\documentclass[11pt,twoside]{article}
\usepackage{latexsym}
\usepackage{amssymb,amsbsy,amsmath,amsfonts,amssymb,amscd}
\usepackage{graphicx,color}
\usepackage{float,url}
\usepackage{cancel}
\usepackage[colorlinks=true]{hyperref}

\newcommand  \ind[1]  {   {1\hspace{-1.2mm}{\rm I}}_{\{#1\} }    }

\newcommand{\commentout}[1]{}
\newcommand{\R}{\mathbb{R}}

\newcommand {\al} {\alpha}
\newcommand {\e}  {\varepsilon}

\newcommand {\Chi} {{\bf \raise 2pt \hbox{$\chi$}} }

\newcommand {\dv}  { {\rm div} }

\newcommand {\calM} { {\mathcal M} }
\newcommand {\f}   {\frac}
\newcommand {\p}   {\partial}
\newcommand{\un}{\underline n}
\newcommand{\dis}{\displaystyle}

\newcommand{\beq}{\begin{equation}}
\newcommand{\eeq}{\end{equation}}
\newcommand{\bea} {\begin{array}{rl}}
\newcommand{\eea} {\end{array}}
\newcommand{\bepa}{\left\{ \begin{array}{l}}
\newcommand{\eepa} {\end{array}\right.}
\newtheorem{theorem}{Theorem}

\newcommand{\qed}{{ \hfill
                       {\unskip\kern 6pt\penalty 500 \raise -2pt\hbox{\vrule\vbox to 6pt{\hrule width 6pt
                       \vfill\hrule}\vrule} \par}   }}
\title{A two species hyperbolic-parabolic model of tissue growth}

\author{Piotr Gwiazda\thanks{Institute of Mathematics, Polish Academy of Sciences,
\'Sniadeckich 8, 00-656 Warszawa, Poland. Email:
pgwiazda@mimuw.edu.pl} \thanks{P.G. and A.\'S.-G. received support from the National Science Centre (Poland),
2015/18/MST1/00075.
}
\and
Beno\^ \i t Perthame\thanks{Sorbonne Universit\'{e}, Universit\'{e} Paris-Diderot SPC, CNRS, INRIA, Laboratoire Jacques-Louis Lions, F-75005 Paris, France. Email : Benoit.Perthame@sorbonne-universite.fr}
\thanks{B.P. has received funding from the European Research Council (ERC) under the European Union's Horizon 2020 research and innovation programme (grant agreement No 740623)}
\thanks{This work has been initiated when this author was visiting the Banach center in Warsaw during Simons Semester partially supported
by the Simons - Foundation grant 346300 and the Polish Government MNiSW 2015-2019
matching fund. }
\and 
Agnieszka \'Swierczewska-Gwiazda\thanks{Institute of Applied Mathematics and Mechanics
University of Warsaw, Banacha 2, 02-097 Warsaw, Poland. Email: aswiercz@mimuw.edu.pl} 
\footnotemark[2]
}
 \date{\today}

\begin{document}
\maketitle
\pagestyle{plain}
%\tableofcontents
\pagenumbering{arabic}

\begin{abstract} 
Models of tissue growth are now well established, in particular in relation to their applications to cancer. They describe the dynamics of cells subject to motion resulting from  a pressure gradient generated by the death  and birth of cells, itself controlled primarily by pressure through contact inhibition. In the compressible regime we consider, when pressure results from the cell densities and when two different populations of cells are considered, a specific difficulty arises from the hyperbolic character of the  equation  for each cell density, and to the parabolic aspect of the equation for the total cell density. For that reason, few a priori estimates are available and discontinuities may occur. Therefore the existence of solutions is a difficult problem. 

Here, we  establish the existence of weak solutions to the model with two cell populations which react similarly to the pressure in terms of their motion but undergo different growth/death rates.  In opposition to the method used in the recent paper \cite{CFSS}, our strategy is to ignore compactness on the cell densities and to prove strong compactness on the pressure gradient. We improve known results in two directions; we obtain new estimates, we treat higher dimension than 1 and we deal with singularities resulting from vacuum.

\end{abstract} 
\vskip .7cm

\noindent{\makebox[1in]\hrulefill}\newline
\section*{Introduction}
%\label{sec:intro}
%-------------------------------------------
%%%%%%%%%%%%%%%%%%%%%%%%%%%%%%%%%%%%%%%%%%%%

The topic of modeling tissue growth has recently progressed with various inputs from physics and mechanics~\cite{PT, BKMP2, GorielyBenAmar2005, RJPJ}.  Models are now used for image-based prediction of cancer growth \cite{BCGR,RSCB}. 
They describe the dynamics of cell number density subject to motion resulting from  a pressure gradient generated by the death  and birth of cells, itself controlled primarily by pressure through contact inhibition. In the compressible regime, pressure results from a combination of the cell densities and controls both the motion through Darcy's law and birth and death of cells according to a finding in \cite{ByDr} and commonly used since then. Models with a single type of cells have been studied recently by many authors, as well as their incompressible limit \cite{PQV, PQTV,  KTu, KPS, MRS2}.  More general formalisms using incompressibility conditions also occur in two phase flows, and they appear, e.g., in oil recovery \cite{ADiB, ALuVi} where each phase has its own pressure. Models may also contain several ``phases'', and have also been widely established and studied \cite{BP,  SCh, CGM, Kpo, BHIM}. For instance,  a specific question is to understand when segregation occurs \cite{ BGH, CFSS}. 
\\

Here, we consider the following compressible two cell population model,  that we state in the full space for the sake of simplicity,
\beq \bepa
\p_t n_1- \dv [n_1 \nabla p] = n_1 F_1(p) +n_2 G_1(p), \qquad x \in \R^d, \; t \geq 0,
\\[10pt]
\p_t n_2- \dv [n_2 \nabla p] = n_1 F_2(p) +n_2 G_2(p),
\eepa
\label{eq1}
\eeq
with
\beq
n:= n_1 + n_2, \qquad p =n^\gamma, \quad \gamma >1.
\label{eq2}
\eeq
 We assume that there is a value $P_H >0$ (the name homeostatic pressure was coined in \cite{RJPJ}) such that the smooth functions $F_i$, $G_i$, describing the division/death rates of cells, satisfy the properties
\beq
F(p):= F_1(p) +F_2(p) \leq 0, \qquad  G(p):= G_1(p) +G_2(p) \leq 0, \qquad \forall p \geq P_H .
\label{as1}
\eeq
We also assume that the initial data $n^0_1$, $n^0_2$, $n^0=  n_1^0 + n_2^0$ satisfy 
\beq
n^0_1 \geq 0, \qquad n^0_2 \geq 0, \qquad p^0:= (n_1^0 + n_2^0)^\gamma \leq P_H ,
\label{as2}
\eeq
\beq
n^0 (1+|x|^2+  |\ln(n^0)| ) \in L^1(\R^d) ,
\label{as3}
\eeq
\beq
\nabla p^0  \in L^2(\R^d) , \qquad  \Delta p^0 \in \calM_{loc}(\R^d), \qquad \big( \Delta p^0 \big)_-  \in L^2_{loc}(\R^d),
\label{as4}
\eeq
where $\calM_{loc}(\R^d)$ refers to the vector space of locally bounded measures. 
At some point, we will also need the restrictions that $\gamma$ is large enough when $d \geq 5$ and that near $p=0$ some cancelation occurs, namely
\beq
\gamma > 2- \f 4 d, \qquad \sup_{0\leq p \leq P_H} \f{ |F(p) - G(p)|^2}{p^{1/\gamma}} \leq C_H .
\label{as5}
\eeq
In words, the total proliferation rates of cells $n_1$ and $n_2$ are the same when $p\approx 0$.

A specific difficulty arises from the hyperbolic character of the  equation  for each cell density $n_i$, and to the parabolic aspect of the total cell density $n$. For example, it is known that solutions $n_1$, $n_2$ may have discontinuities.  For that reason, the existence of solutions is a difficult problem by lack of strong a priori estimates. Also we cannot hope for strong solutions in general.  
Here, we  establish the existence of weak solutions.  In opposition to the method used in the recent paper \cite{CFSS}, our strategy is to ignore compactness on the cell densities and to prove strong compactness on the pressure gradient. Therefore, we improve known results in two directions; we treat higher dimension than 1 as in \cite{CFSS} and we deal with vacuum while \cite{BHIM} only considers uniformly positive and smooth solutions.

%----------------------------------
\begin{theorem}[A priori estimates]
With the assumptions \eqref{as1}--\eqref{as4},  the following estimates hold true for all $T>0$ with constants $C(T)$ which only depend on the bounds in the above assumptions
\beq
n(x,t) \geq 0 , \qquad \int_{\R^d} n(x,t)  \,dx \leq Ce^{Ct}, \qquad  p(x,t) \leq P_H,
\label{est1}
\eeq
\beq
\int_0^T \int_{\R^d}  \f{| \nabla p |^2}{p^{1-1/\gamma}}  \,dx dt \leq C(T).
\label{est2}
\eeq
Assuming also \eqref{as5}, we have for $t\in (0,T)$,  and with  $\phi(\cdot) \in C^2_{\rm Comp} (\R^d)$ a localizing function,
 \beq 
\dis \int_{\R^d}   |\Delta p (t) |_-^2 \phi(x)  \,dx \leq C(T), \qquad  \int_0^T  \int_{\R^d}   |\Delta p (t) |_-^3 \phi(x)  \,dx dt \leq C(T),
\label{est5}
\eeq
 \beq 
\int_{\R^d}   |\Delta p (x,t) | \;  \phi(x)   \,dx \leq C(T) .
\label{est6}
\eeq
\label{th:estimates}
\end{theorem}
%-------------------------------
Since our framework includes the Barenblatt solutions, see \cite{Vbook}, we know that these estimates are sharp in the sense that $\Delta p$ may be a singular measure supported by the free boundary. Note that $p$ being bounded, the estimate~\eqref{est2} also gives an $L^2_{t,x}$ bound on $\nabla p$.  Another a priori estimate is also available, which we do not use in the subsequent results, and that we postpone to the Appendix.
\\

As a consequence of the estimates in Theorem~\ref{th:estimates}, we establish the following stability result
%
%---------------------------------------------------
\begin{theorem} [Stability of weak solutions] Assume \eqref{as1} and \eqref{as5} and 
that the family of initial data satisfies, with uniform bounds, the assumptions \eqref{as2}--\eqref{as4}. Then, the corresponding weak solutions $n_i^\e$, with the above bounds true, satisfy after extraction of subsequences, 
$$
n_i^\e  \rightharpoonup n_i, \quad \hbox{ in } \; L^\infty \big((0,T)\times \R^d \big)-w*, \qquad i=1, \, 2,
$$
$$
n^\e \to n, \qquad p^\e \to p , \quad \hbox{ in } \; L^q \big((0,T)\times \R^d \big), \; 1\leq q < \infty,
$$
$$
\nabla p^\e \to \nabla p  , \quad \hbox{ in } \; L^2 \big((0,T)\times \R^d \big)
$$
and $n_1$, $n_2$, $p$ satisfy, in the weak sense, the system \eqref{eq1}--\eqref{eq2} with initial data $n^0_1$, $n^0_2$.
\label{th:stability}
\end{theorem} 
%---------------------------------------------------
Finally, these two results lead us to the existence theorem, which is the main result of the current paper.
\begin{theorem}[Existence of weak solutions]
With the assumptions  of Theorem \ref{th:stability},  there exists a weak solution $n_1,n_2,p\in L^\infty((0,T)\times\R^d)$
to the system \eqref{eq1}--\eqref{eq2}, i.e., for $i=1,\, 2$
\beq
\int_0^T\int_{\R^d} \left[-n_i\p_t\psi+n_i\nabla p.\nabla \psi- \big(n_1F_i(p)+n_2G_i(p)\big) \psi \right]  \,dxdt =\int_{\R^d}n_i^0 \psi(0) \,dx
\eeq
holds for all $\psi\in C^1_{\rm Comp}(\R\times\R^d)$ and relations \eqref{eq2} hold a.e. in $(0,T)\times\R^d$.
\label{th:existence}
\end{theorem}

\medskip

The main observation is that, while our problem is of hyperbolic nature,  we can take advantage of informations coming from the parabolic equation on $n$
\beq 
\p_t n - \dv [n \nabla p] = n_1 F(p) +n_2 G(p)=: n R(c_1, c_2, p) , \qquad x \in \R^d, \; t \geq 0,
\label{eqn}
\eeq
where, following \cite{CFSS}, we define 
$$
c_i= \f{n_i}{n}\leq 1  \quad \hbox{and } \;  c_i (x,t)=0 \; \hbox{ when } \; n(x,t)=0,
$$ 
\beq 
R = c_1 F(p) +  c_2 G(p) \in L^\infty . 
\label{defR}
\eeq

Next, multiplying equation~\eqref{eqn} with $p'(n)$,  we compute that $p$ satisfies
\beq
\p_t p - |\nabla p|^2 - \gamma p \Delta p = \gamma p R.
\label{eqp}
\eeq

It is also useful for later purpose to state the equation for the $c_i$'s
\beq
\p_t  c_i - \nabla p.\nabla c_i = c_1F_i(p)+c_2G_i(p) - c_iR.
\label{eqc}
\eeq
To obtain the equation for $c_i$ we multiply the equation for $n_i$ with $1/n$ and add it to the equation for $n$ multiplied by $-n_i/n^2$.  Indeed, observe that
\[
-\frac{1}{n}\dv[n_i\nabla p]+\f{n_i}{n^2}\dv[n\nabla p]=-\f{1}{n}\nabla n_i .\nabla p-\f{n_i}{n}\Delta p+\f{n_i}{n^2}\nabla n.\nabla p+\f{n_i}{n^2}n\Delta p=-\nabla(\f{1}{n}n_i).\nabla p
\]
and the remaining terms are immediate. 
\\

The rest of  the paper is devoted to the proofs of these three theorems which we perform in the three next sections. Some  remarks and open problems are commented in the conclusion.

%-----------------------------------------------------------
\section{Proof of Theorem \ref{th:estimates}}
%----------------------------------------------------------
%%%%%%%%%%%%%%%%%%%%%%%

The first estimates come from the balance law expressed by the equation  \eqref{eq1} and from the maximum principle for equation \eqref{eqp}.  One easily gets, by integrating~\eqref{eqn} over $\R^d$ and using the Gronwall inequality, that
\beq
\int_{\R^d}n(t,x)\,dx\le \int_{\R^d}n^0(x) \,dx \exp(t\|R\|_\infty).
\eeq
To show the uniform bound on $p$ we multiply~\eqref{eqp} with $(p-P_H)_+$. Observe that for any $\eta\in C^2$ it holds $\Delta \eta(p)=\Delta p\eta'(p)+\eta''(p)|\nabla p|^2$, which allows us to handle the highest order term with $\eta'(p)=p(p-P_H)_+$. Thus we get
\beq
\f{1}{2}\partial_t(p-P_H)_+^2+[\gamma\eta''(p)-(p-P_H)_+] |\nabla p|^2-\gamma\Delta\eta(p)=\gamma p(p-P_H)_+R(p).
\eeq
We integrate over $\R^d$ and observe that as $\gamma>1$, thus $\gamma\eta''(p)-(p-P_H)_+\ge0$ and
\beq
\f{1}{2}\f{d}{dt}\int_{\R^d}(p-P_H)_+^2 \,dx\le\gamma\int_{\R^d}p(p-P_H)_+R(p)\le0,
\eeq
where the last inequality follows from \eqref{as1} and \eqref{as2}. 
\\

Before passing to the next estimate, on $|\nabla p|^2$,  let us observe that the second moment of $n$ is bounded. 
Indeed, multiplying~\eqref{eqn} with $x^2\Phi_L$, where $\Phi_L(|x|)$ is a radially symmetric smooth function, which vanishes outside the ball of radius $L+1$, equals to $1$ on the ball of radius $L$, with $\nabla \Phi_L$ and
$\Delta\Phi_L$ bounded uniformly in $L$;  and integrating by parts over $\R^d$ gives 
\beq
\f{d}{dt}\int_{\R^d}x^2n\Phi_L  \,dx+\f{\gamma}{\gamma+1} \int_{\R^d}n^{\gamma+1}(x^2\Delta \Phi_L+4x\nabla \Phi_L+2\Phi_L)  \,dx=\int_{\R^d}x^2nR \Phi_L \,dx.
\eeq
Since $n^\gamma\le P_H$, we can furthermore obtain
\beq
\f{d}{dt}\int_{\R^d}x^2n  \Phi_L\,dx\le C\int_{\R^d}x^2n \Phi_L \,dx+\int_{\R^d}n\Phi_L \,dx
+\int_{\{L\le|x|\le L+1\}}nx |\nabla \Phi_L|\,dx+\int_{\{L\le|x|\le L+1\}}n x^2|\Delta\Phi_L|\,dx,
\eeq
where the constant $C$ depends on 
$\|R\|_\infty$, $\gamma$ and $P_H$.
As $n\in L^1$, we can claim that the  second term on the right-hand side is bounded. 
For the moment let us assume that $n$ vanishes sufficiently fast at infinity, what we will prove later.
Then, by the  Lebesgue dominated convergence theorem we can pass to the limit in terms containing $\Phi_L$ and using that $n$ vanishes for large $x$ and both $\nabla \Phi_L$ and $\Delta\Phi_L$ are uniformly bounded, we show that the last two terms on the right-hand side vanish. 
We complete the estimate by applying the Gronwall inequality. 
\\

The estimate \eqref{est2} comes from  the entropy relation. We multiply~\eqref{eqn} with $\Phi_L \ln n$, where $\Phi_L$ is the same truncation function as above, integrate over $\R^d$, and find
\begin{equation}
\begin{aligned}
\frac{d}{dt}& \int_{\R^d} n( \ln (n) -1) \Phi_L\,dx+  \f{1}{\gamma}\int_{\R^d} p^{-1+1/\gamma}  |\nabla p|^2 \Phi_L\,dx\\
&-\f{1}{\gamma-1}\int_{\{L\le |x|\le L+1\}}p^{1-\frac{1}{\gamma}}(\ln (p)-1)\Delta\Phi_L\,dx
=   \int_{\R^d} n \ln(n) \;  R\Phi_L\,dx.
\end{aligned}
\end{equation}
It is easy to observe that if the function $n$, and thus also $p$ vanishes sufficiently fast, then again the integral over the annulus vanishes as $L\to\infty$ and we conclude~\eqref{est2}. For that purpose we recall also that a control of the second moment in $x$ is used here to control the negative values of $n \ln(n)$. Indeed, observe that
\beq
\int_{\R^d}n|\ln (n)|  \,dx=\int_{\R^d} n\ln (n)  \,dx -\int_{\{x\in\R^d: n<1\}}n\ln (n)  \,dx\le \int_{\R^d} n\ln (n)  \,dx+ 2\int_{\R^d}n|x|^2  \,dx+c .
\eeq
And the above inequality allows us to get 
\beq
\int_{\R^d}n|\ln (n)|  \,dx\le \|R\|_\infty\int_0^T\int_{\R^d}n|\ln (n)| +2\int_{\R^d}n|x|^2 +\int_{\R^d}n \,dx+\int_{\R^d}n^0(\ln (n^0)-n^0) \,dx+c.
\eeq
We complete the estimate \eqref{est2}  using the Gronwall lemma. 
\\

Finally, the fundamental estimates \eqref{est5} come from Aronson and Benilan's method~\cite{AB, Vbook} for the porous media equation with several  adaptations. Firstly, and this s a new feature here, we weaken their estimate to $L^2$  rather than $L^\infty$. Secondly, we need to localize the estimate in space. Thirdly, we adapt the functional under consideration using also the idea from \cite{PQV},  and we do not work  directly with  $\Delta p$ but with 
$$
w = \Delta p + R, \qquad \p_t p = |\nabla p |^2 + \gamma p w.
$$
We compute 
$$
\p_t \Delta p = 2 (\p_{ij} p)^2  + 2 \nabla p \nabla \Delta p + \gamma \Delta (pw),
$$
\begin{equation}
\begin{aligned}
\p_t  R &= R_{c_1} \p_t c_1 +  R_{c_2 }\p_t c_2 +  R_p \p_t p
\\
&= F(p) [\nabla c_1 . \nabla p + \xi_1]  + G(p)  [\nabla c_2 . \nabla p + \xi_2] + R_p  [ |\nabla p |^2 + \gamma p w ]
\end{aligned}
\label{Rt}
\end{equation}
where 
$$
\xi_i:=c_1F_i(p)+c_2G_i(p)+c_iR, \qquad i=1,\, 2,
 $$ 
are the right-hand sides from  equations \eqref{eqc}.
Therefore, since $ c_1 +  c_2=1$, we find
$$
\p_t w =  2 (\p_{ij} p)^2 +  2 \nabla p \nabla \Delta p + \gamma \Delta (pw) + [ F(p)  -G(p) ] \nabla c_1 . \nabla p + R_p  [ |\nabla p |^2 + \gamma p w ] +Bdd_1
$$
with ``$Bdd$'' terms which are bounded in $L^\infty$, 
$$
Bdd_1:=F(p)\xi_1+G(p)\xi_2,
$$
but this may change from line to line. 
Since
$$
\nabla p .\nabla\Delta p=\nabla p.\nabla (w-R)  =   \nabla p.\nabla w-\dv (R\nabla p)+R(w-R)
$$
 this is also
\begin{equation}\begin{aligned}
\p_t w &\geq   \f 2 d (\Delta p)^2 +  2 \nabla p \nabla w - 2 \dv (  R \nabla p ) + 2 R (w-R)+\gamma \Delta (pw)\\
& + [ F(p)  -G(p) ] \nabla c_1 . \nabla p + R_p  [ |\nabla p |^2 + \gamma p w ] +Bdd_1.
\end{aligned}\end{equation}
 The negative part, that we denote by  $| w |_-$,  therefore satisfies, with $sgn_-  := \ind{w<0}$,
\begin{equation}\begin{aligned}
\p_t | w |_- \leq & -\f 2 d |w|_-^2 +  2 \nabla p \nabla | w|_-  + 2 sgn_- \dv (  R \nabla p ) + 2(1-\f 2d)  R |w |_- +  \gamma\Delta (p |w|_- ) 
\\ 
& - sgn_- [ F(p)  -G(p) ] \nabla c_1 . \nabla p - sgn_-  R_p   |\nabla p |^2 + \gamma R_p p |w |_-  +Bdd
\end{aligned}\end{equation}
using that $-\f 2d(w-R)^2=-\f 2d w^2+\f 4dRw-\f 2dR^2$, where the last term we include within bounded terms that we still gather in~ $Bdd$.

We reorganize this inequality as (here the parameter $\al >0$ can be chosen as small as we wish)
\[\begin{aligned}
\p_t | w |_- \leq & - (\f 2 d -\al) |w|_-^2 +  2 \nabla p \nabla | w|_-   +  \gamma\Delta (p |w|_- ) + 2 sgn_- \dv (  R \nabla p )
\\ 
& - sgn_- [ F(p)  -G(p) ] \nabla c_1 . \nabla p - sgn_-  R_p   |\nabla p |^2   +Bdd \; . 
\end{aligned}\]
Notice that, above, we have applied the Young inequality to the terms  $2(1-\f 2d)  R |w |_-$ and $\gamma R_p p |w |_-$, that is 
$$
2R|w|_-(1-\f 2d)\le \f{\alpha}{2}|w|_-^2+c(\alpha)R^2, \qquad 
R_p\gamma p |w|_-\le\f{\alpha}{2}|w|_-^2+c(\alpha)|R_p\gamma p|^2,
$$ 
and thus the term $Bdd$ is given by
$$
Bdd:=F(p)\xi_1+G(p)\xi_2-\f 2dR^2+c(\alpha)(R^2+|\gamma R_p p|^2).
$$

We need to localize and use a nonnegative, compactly supported, smooth test function $\Phi$ to compute 
\beq
\f{d}{dt} \int_{\R^d} \f{| w |_-^2 }{2} \Phi \, \,dx \leq I + II .
\label{eq:wsquare}
\eeq
The terms $II$ are those with $\nabla \Phi$ which are certainly better because they contain one less derivative in the unknowns. 
\\

For the difficult term we have, after several integrations by parts,   in particular to eliminate derivatives in $c_1$ which are the worse,  
\begin{equation}\begin{aligned}
I \leq   & -(\f 2 d - \al)   \int_{\R^d}  \Phi  |w|_-^3 \,dx -  \int_{\R^d}   \Phi | w|_-^2  \Delta p  \,dx  - \gamma  \int_{\R^d}   \Phi  p | \nabla  |w|_- |^2 \,dx + \f \gamma 2  \int_{\R^d}  \Phi  | w|_-^2  \Delta p  \,dx
\\ 
& - 2  \int_{\R^d}  \Phi R \nabla  | w|_-  .  \nabla p \,dx + \int_{\R^d} \Phi c_1    [ F(p)  -G(p) ] | w|_-  \Delta p \,dx
\\ 
& + \int_{\R^d}  \Phi   c_1   [ F(p)  -G(p) ]  \nabla | w|_-  . \nabla p \,dx+ C \int_{\R^d}  \Phi  | w|_-   |\nabla p |^2 \,dx +\int_{\R^d}  \Phi  | w|_- Bdd \,dx\; .
\end{aligned}\end{equation}
where $C$ is constant which here takes into account $ c_1 (F'-G')$ as well as $R_p$
and which is changing from line to line below. 

The linear and quadratic terms in $|w|_-$ are not a problem because the dominant term contains a cubic power of $|w|_-$. 
We observe about the third-power terms that 
$$
\int_{\R^d}\Phi|w|_-^2\Delta p\,dx=\int_{\R^d}\Phi|w|_-^2(w-R)\,dx=-\int_{\R^d}\Phi|w|_-^3\,dx-\int_{\R^d}\Phi|w|_-^2R\,dx.
$$

Because $R = G(p) +c_1   [ F(p)  -G(p) ] $, we also have
\begin{equation}\begin{aligned}
I \leq  & -( \f \gamma 2 + \f 2  d -1-\al)   \int_{\R^d}   \Phi  |w|_-^3\,dx     - \gamma  \int_{\R^d}   \Phi  p | \nabla  |w|_- |^2 \,dx
\\ 
&- 2  \int_{\R^d}  \Phi  G  \nabla  | w|_-   \nabla p\,dx + \int_{\R^d} \Phi c_1    [ F(p)  -G(p) ] | w|_-  \Delta p \,dx
\\ 
& - \int_{\R^d}   \Phi   c_1 [ F(p)  -G(p) ]  \nabla | w|_- . \nabla p \,dx+ C \int_{\R^d}  \Phi  | w|_-   |\nabla p |^2 \,dx
\\ 
&+(1-\f{\gamma}{2})\int_{\R^d}\Phi|w|_-^2R\,dx+\int_{\R^d} \Phi | w|_- Bdd\,dx \; ,
\end{aligned}\end{equation}
it holds
$$
\int_{\R^d} \Phi c_1    [ F(p)  -G(p) ] | w|_-  \Delta p\,dx=-\int_{\R^d}\Phi c_1[F(p)-G(p)]|w|_-^2\,dx-\int_{\R^d}\Phi c_1[F(p)-G(p)]|w|_-R\,dx
$$
and estimating further (the second term of the above we include already in $\int_{\R^d} \Phi|w|_-Bdd\,dx$),
\begin{equation}\begin{aligned}
I \leq  & -( \f \gamma 2 + \f 2  d -1-\al)   \int_{\R^d}   \Phi  |w|_-^3\,dx     - \gamma  \int_{\R^d}   \Phi  p | \nabla  |w|_- |^2 \,dx
\\ 
&+ 2  \int_{\R^d}  \Phi  G  | w|_-   \Delta p\,dx+2\int_{\R^d} \Phi G'(p)|w|_-|\nabla p|^2\,dx
\\ 
&+  \f 12 \int_{\R^d}   \Phi  p c_1^2 | \nabla  |w|_- |^2 \,dx +\f12  \int_{\R^d}   \Phi  \f{[ F(p)  -G(p) ]^2}{p} |\nabla p|^2 \,dx+ C \int_{\R^d}  \Phi  | w|_-   |\nabla p |^2\,dx  \\ 
&+(1-\f{\gamma}{2})\int_{\R^d}\Phi|w|_-^2R\,dx-\int_{\R^d}\Phi c_1[F(p)-G(p)]|w|_-^2\,dx+\int_{\R^d} \Phi | w|_- Bdd\,dx
\;.
\end{aligned}\end{equation}
We arrive at the final form, using the constant $C_H$ in~\eqref{as5}, 
\begin{equation}\begin{aligned}
I \leq  & -( \f \gamma 2 + \f 2  d -1-\al)   \int_{\R^d}   \Phi  |w|_-^3 \,dx    - (\gamma -\f 12)  \int_{\R^d}   \Phi  p | \nabla  |w|_- |^2  \,dx
\\ 
&+C\int_{\R^d}\Phi|w|_-^2\,dx
 +\f{C_H}{2}  \int_{\R^d}   \Phi  \f{ |\nabla p|^2}{p^{1-1/\gamma} } \,dx+ C \int_{\R^d}  \Phi  | w|_-   |\nabla p |^2\,dx +\int_{\R^d}  \Phi | w|_- Bdd\,dx \; .
\end{aligned} \end{equation}
where a constant $C$ standing next to the integral $ \int_{\R^d}  \Phi  | w|_-   |\nabla p |^2\,dx $ takes into account also $2G'$. 
The first two terms on the right-hand side are the ``good terms''. Indeed, they have a good sign, unless $\gamma>2(1-\f 2d)$, which is automatically satisfied as we assume $\gamma>1$ in $d\le 4$ and in higher dimensions implies higher requirement of the exponent $\gamma$ as stated in~\eqref{as5}. 
The difficult term is at the highest order 
\beq
\int_{\R^d}  \Phi  | w|_-   |\nabla p |^2\,dx = -  \int_{\R^d}  \Phi  [ \nabla | w|_-    p \nabla p + | w|_-   p \Delta p ]\,dx.
\label{highest}
\eeq
which is under control. Indeed, 
$$
 \int_{\R^d}  \Phi   \nabla | w|_-    p \nabla p\,dx\le (\gamma-\f{1}{2})\int_{\R^d} \Phi p |\nabla|w|_-|^2\,dx+c(\gamma)\int_{\R^d}\Phi p|\nabla p|^2\,dx.
$$
Here, the first term on the right hand side just cancels the second ``good term''and the second is bounded.

 For the second term of the right-hand side of~\eqref{highest} we have
$$
 \int_{\R^d}  \Phi  | w|_-   p \Delta p \,dx= \int_{\R^d}  \Phi  | w|_-   p (w-R)\,dx= -\int_{\R^d}  \Phi  | w|_-^2   p \,dx+ \int_{\R^d}  \Phi  | w|_-   p R\,dx
$$
where both these terms are under control to give the final estimate
\beq
\begin{aligned}
I \leq  -( \f \gamma 2 + \f 2  d -1-\al)   \int_{\R^d}   \Phi  |w|_-^3\,dx  + C. 
\end{aligned}
\label{eq:wcubic} \eeq

The terms containing gradient of $\Phi$ are collected in $II$
\begin{equation}\begin{aligned}
II&= \int_{\R^d}\nabla p|w|_-^2\nabla \Phi\,dx-\gamma\int_{\R^d}\nabla p|w|_-^2\nabla \Phi\,dx
-\gamma\int_{\R^d} p\nabla(|w|_-^2)\nabla \Phi\,dx
-2\int_{\R^d} R\nabla p |w|_-\nabla \Phi\,dx\\ 
&+\int_{\R^d}[F(p)-G(p)]c_1\nabla p|w|_-\nabla \Phi\,dx+2\int_{\R^d}\nabla\Phi G|w|_-\nabla p\,dx-\int_{\R^d}\nabla \Phi |w|_-p\nabla p\,dx
\end{aligned}\end{equation}
and they  all do not bring additional difficulties. 

Therefore, using the inequality \eqref{eq:wsquare} and the negative sign in the right hand side of \eqref{eq:wcubic}, we obtain the a priori estimates  announced in~\eqref{est5}.  The $L^1$ bound for $\Delta p$ in ~\eqref{est6} is a simple consequence because
$\int \Phi \Delta p dx $ is bounded, therefore $\int \Phi |\Delta p|_+ dx $ is controlled by $\int \Phi |\Delta p|_- dx $ which itself is controlled thanks to~\eqref{est5}.  

%----------------------------------------------------------
%%%%%%%%%%%%%%%%%%%%%%%
\section{Proof of Theorem~\ref{th:stability}}
\label{sec:stability}
%----------------------------------------------------------
%%%%%%%%%%%%%%%%%%%%%%%

The goal here is to explain the main compactness argument which is used to pass to the limit in an approximate sequence. As in \cite{CFSS}, the compactness in time is a major issue.
\\

Weak convergence of the quantities $n_i^\e$ follows from the bound in $L^\infty$. The strong convergence of $p^\e$ follows from compactness by Sobolev injections. Indeed, on the one hand,  we control $\int_0^T \int_{\R^d} |\nabla p^\e|^2 dx dt$ from \eqref{est2} because the pressure is bounded by $P_H$. On the other hand, we may  win time compactness by the Lions-Aubin Lemma using  equation~\eqref{eqp}, which also reads
\beq
\p_t p^\e =  (\gamma -1 )  |\nabla p^\e |^2 +  \f{\gamma}{2}  \Delta p^{\e2 }+ \gamma p^\e R^\e ,
\label{eqpdiv}
\eeq
and the space compactness on $p^\e$ together with the known bounds provide time compactness. Therefore  the expression $n^\e =(p^\e)^{1/\gamma}$ shows that we may also extract a sub-sequence of $n^\e$ which converges.
\\

The strong compactness for $\nabla p^\e$ is more involved. It mainly relies on the second estimate ~\eqref{est5} which provides the space compactness, still by the Sobolev embedding theorems. Indeed, the control in $L^1_{\rm loc}$ of $\Delta p^\e$ is enough for compactness of $\nabla p^\e$, a fact which can be inferred from the representation formula for the solution of the Laplace equation.

For time compactness of $\nabla p^\e$, we write, using \eqref{eqpdiv},
$$
\p_t \nabla p^\e =  \nabla \big[ (\gamma -1) |\nabla p^\e |^2 + \gamma p^\e R^\e \big]+  \f{\gamma}{2}  \nabla \Delta p^{\e2 } .
$$
Again, we  know the local space compactness of $ \nabla p^\e$ from the previous paragraph, the right hand side is a sum of space derivatives of bounded functions, therefore we may apply the Lions-Aubin compactness argument and find that $\nabla p^\e$ is compact in space and time.
\\

To pass to the limit in the equations is now easy. All the nonlinear terms, that are
$$
n^\e_i \nabla p^\e, \qquad n^\e_i F_j(p^\e), \qquad n^\e_i G_j(p^\e), 
$$
 have limits as products of weak limits of $n^\e_i$ by strong limits of $p^\e$ and $\nabla p^\e$. This completes the proof of Theorem~\ref{th:stability}.
 \qed

\medskip

At this stage, let us point out that our strategy  differs deeply from that in \cite{CFSS} based on BV estimates for the quantities $c^\e_i$ in one dimension. This estimate is somehow sharp since examples with discontinuities on the $n^\e_i$ are known. Also the method for time compactness is very different since \cite{CFSS}  use a control of the Wasserstein distance.

%----------------------------------------------------------
%%%%%%%%%%%%%%%%%%%%%%%
\section{Proof of Theorem~\ref{th:existence}}
%----------------------------------------------------------
%%%%%%%%%%%%%%%%%%%%%%%

We already have a priori estimates and a  weak sequential stability result, thus to complete the existence proof we need to construct an approximate system compatible with these estimates. We do that in two steps. Firstly, we make positive the initial data and prove a control from below by a (small) Gaussian. Secondly, we introduce a uniform parabolic regularization.
\\

\noindent {\em First step. A regularized problem with a positive control from below.}  We show that the function 
\beq
\underline{n}(t,x)=\underline c\exp\left({-\f{|x|^2}{2}-ct}\right)
\label{subsol}
\eeq
is a subsolution to equation~\eqref{eqn} if we choose $c$ sufficiently large.
Since $\partial_t \underline n=-c\underline n$, $\nabla \underline n=- x \underline n $ and  $\Delta \underline n =d \underline n + |x|^2 \underline n$, we may insert~\eqref{subsol} into the equation for $n$ and as we search for a subsolution, we change  equality to inequality. We obtain
\beq
-c\un-\gamma(\gamma +1)\un^{\gamma+1}|x|^2+\gamma \un^{\gamma+1}-R\un\le 0
\eeq
which holds true  choosing $c$ large enough so that the inequality is satisfied
\beq
-c-\gamma(\gamma +1)\left[\exp\left({-\f{|x|^2}{2}-ct}\right)\right]^{\gamma}|x|^2+\gamma \left[\exp\left({-\f{|x|^2}{2}-ct}\right)\right]^{\gamma}+\|R\|_\infty\le 0.
\eeq
It is now a matter of standard estimates, \cite{Vbook}, to obtain that if we start with the specific initial condition larger than $n^0$, we will call it
$n_\delta^0:=n^0+\delta\exp(\f{-|x|^2}{2})$, with $\delta>0$, then the solution to the problem will be larger than the  subsolution $\un$ given by \eqref{subsol} with $c$ large enough.

Thus our first  approximation step is to replace an initial data $n^0$ by $n_\delta^0$ as it is introduced above. In a consequence the corresponding solution $n_\delta$, as well as $p_\delta$ are locally bounded away from zero, we call these bounds $\un_\delta$ and $\underline p_\delta$. The solution
has a regularity $L^q(0,T;W^{2,q}_{loc}(\R^d))$, see~\cite{Vbook} for details and the method in \cite{BHIM}.

Note that an analogue maximum estimate can be proven to provide a bound from above and justify that $n$ vanishes at infinity, what we announced earlier.  
\\

\noindent  {\em Second step. A uniformly parabolic approximation.} We consider  the system of equations, which consists of a parabolic equation for    $n$ and hyperbolic equations for  $c_1$ and $c_2$.
Thus we construct a parabolic approximation of the equation for $c_i$, $i=1,2$. Let $\e>0$, 
\beq
\p_t c_i^\e-\nabla p^\e . \nabla c_i^\e-\e \,\dv [p^\e\nabla c_i^\e]=c_1^\e F_i(p^\e)+c_2^\e G_i(p^\e) - c_i^\e R(p^\e).
\label{eqceps}
\eeq
Note that all the quantities are for simplicity labelled only with $\e$, but they depend both on $\e$ and $\delta$, i.e. $p^\e:=p^{\e,\delta}$ as well as the other quantities. 
We proceed now as follows: We solve a parabolic system consisting of \eqref{eqn} and \eqref{eqceps} with initial data $p^0_\delta$ and $\frac{n_{i,\delta}^0}{n^0_\delta}$, completed with the relation $p^\e=(n^\e)^\gamma$. 
The equations \eqref{eqceps} allow to observe the crucial property, which $c_i$ possessed and  which was the only information on these quantities used in a priori estimates. Indeed, adding the equations on $c_i^\e$, we keep the fundamental relationship $c_1^\e+c_2^\e \equiv 1$ thanks to the definition of $R$ in \eqref{defR}, since initially $c_1^\e+c_2^\e = 1$. Finally, with the fully parabolic framework at hand, it is in the folklore of the domain  to obtain the existence of the coupled problem between $n^\e $ and $c_i^\e$.

Next, we notice  that all the  a priori bounds used to pass to the limit are true. Multiplying with $c_i^\e|w^\e|_-\Phi$ and integrating over $(0,T)\times\R^d$ gives
\begin{equation}
\begin{aligned}
\int_{\R^d}|c_i^\e|^2|w^\e|_- \Phi  \,dx&+\int_0^T\int_{\R^d}\nabla p^\e . \nabla c_i^\e c_i^\e |w^\e|_-\Phi  \,dxdt\\
&+\e\left[\int_0^T\int_{\R^d}p^\e  |\nabla c_i^\e|^2 |w^\e|_-\Phi  \,dxdt+\int_0^T\int_{\R^d}p^\e  \nabla c_i^\e \nabla |w^\e|_-\Phi  \,dxdt\right]\\
&=\int_0^T\int_{\R^d} [c_1^\e F_i(p^\e)+c_2^\e G_i(p^\e)+c_i^\e R(p^\e)]|w^\e|_-c_i^\e \Phi  \,dxdt+\int_{\R^d}|c_i^\e(0)|^2 |w^\e|_-\Phi  \,dx.
\end{aligned}
\end{equation}

%We estimate 
%$$
%|\int_0^T\int_{\R^d}\nabla p^\e . \nabla c_i^\e c_i^\e \Phi  \,dxdt|
%\le \|c_i^\e\|_\infty\left(\f{2\|c_i^\e\|_\infty}{2\e\underline{p}_\delta}\int_0^T\int_{\R^d}|\nabla p^\e|^2\Phi + \f{\e\underline{p}_\delta}{2\|c_i^\e\|_\infty}\int_0^T\int_{\R^d}|\nabla c_i^\e|^2 \Phi  \,dxdt\right)
%$$
%Since $p(x,t)\ge \underline{p}_\delta$ thus
%\beq
%\int_{\R^d}|c_i^\e|^2 \Phi  \,dx+
%+\f{\e\underline{p}_\delta}{2}\int_0^T\int_{\R^d} |\nabla c_i^\e|^2 \Phi  \,dxdt\le \f{C}{2\e\underline{p}_\delta}\int_0^T\int_{\R^d}|\nabla p^\e|^2\Phi   \,dxdt
%+\int_{\R^d}|c_i^\e(0)|^2 \Phi  \,dx+C
%\eeq
We observe that 
\begin{equation}\begin{aligned}
\int_0^T\int_{\R^d}\nabla p^\e . \nabla c_i^\e c_i^\e|w^\e|_- \Phi  \,dxdt&=\f 12\int_0^T\int_{\R^d}\nabla p^\e . \nabla |c_i^\e|^2|w^\e|_- \Phi  \,dxdt\\ 
&=-\f 12\int_0^T\int_{\R^d}\Delta p^\e  |c_i^\e|^2 \Phi  \,dxdt-\f 12\int_0^T\int_{\R^d}\nabla p^\e  |c_i^\e|^2 \nabla \Phi  \,dxdt\\ 
&\quad-\f 12\int_0^T\int_{\R^d}\nabla p^\e  |c_i^\e|^2 \nabla|w^\e|_- \Phi  \,dxdt.
\end{aligned}\end{equation}
Consequently, we find that
\begin{equation}\begin{aligned}
\e\left[\int_0^T\int_{\R^d}\right.&\left.p^\e  |\nabla c_i^\e|^2 |w^\e|_-\Phi  \,dxdt+\int_0^T\int_{\R^d}p^\e  \nabla c_i^\e \nabla |w^\e|_-\Phi  \,dxdt\right]\\ 
&\le C\int_0^T\int_{\R^d}|\Delta p^\e|   \Phi  \,dxdt
+C\int_0^T\int_{\R^d}|\nabla p^\e | \nabla \Phi  \,dxdt+C\int_0^T\int_{\R^d}\nabla p^\e  |c_i^\e|^2 \nabla|w^\e|_- \Phi  \,dxdt\\ 
&+C\int_0^T\int_{\R^d} |w^\e|_-c_i^\e \Phi  \,dxdt+\int_{\R^d}|c_i^\e(0)|^2 |w^\e|_-\Phi  \,dx.
\end{aligned}\end{equation}
Therefore the first integral on the right-hand side can be again 
 estimated by
$$
\int_0^T\int_{\R^d}|\Delta p^\e|\Phi\, dxdt\le \int_0^T\int_{\R^d}(|w^\e|_-+R)\Phi\, dxdt.
$$

This approximation step will affect our a priori estimates on the level of using computations~\eqref{Rt}, as additional terms related with parabolic approximation, which are estimated above, will appear. Taking these into account, we may pass to the limit as in Section~\ref{sec:stability}, first with $\e\to 0$  and no major difficulty arises. Thus we obtain a limit system for $n$ and $c_i$'s, but still we lack the information whether the  equations  for $n_i$ are satisfied in distributional sense. To recover this  we multiply the equations for $c_i^\e$ with $n^\e$ and add the equation for $n^\e$ multiplied with $c_i^\e $. Let us then define $n_i^\e:=c_i^\e n^\e$ and observe that this operation will lead us to equations for $n_i^\e$ 
\beq 
\p_t n_i^\e- \dv [n_i^\e \nabla p^e]-\e  \dv[p^\e\nabla c_i^\e]n^\e= n_1^\e F_i(p^\e) +n_2^\e G_i(p^\e)
\eeq
The only term that needs to be discussed is $\e \dv[p^\e\nabla c_i^\e]n^\e$. To show that this term vanishes in a limit observe that
\begin{equation}\begin{aligned}
\int_{\R^d}\dv[p^\e\nabla c_i^\e]n^\e \phi dx&=-\int_{\R^d} p^\e\nabla c_i^\e.\nabla [(p^\e)^\frac{1}{\gamma}] \phi dx-
\int_{\R^d} p^\e\nabla c_i^\e (p^\e)^\frac{1}{\gamma} \nabla \phi dx\\
&=-\frac{1}{\gamma}\int_{\R^d}(p^\e)^\frac{1}{\gamma}\nabla p^\e. \nabla c_i^\e \phi dx- \int_{\R^d}  (p^\e)^{\frac{1}{\gamma}+1}\nabla c_i^\e.\nabla \phi dx\\
&=\frac{1}{\gamma}\int_{\R^d}(p^\e)^\frac{1}{\gamma}\Delta p^\e  c_i^\e \phi dx+ 
\frac{1}{\gamma}\int_{\R^d}\nabla ((p^\e)^\frac{1}{\gamma}).\nabla p^\e  c_i^\e \phi dx+
\frac{1}{\gamma}\int_{\R^d}(p^\e)^\frac{1}{\gamma}\nabla p^\e  c_i^\e \nabla \phi dx\\
&+
\int_{\R^d}  (p^\e)^{\frac{1}{\gamma}+1} c_i^\e\Delta \phi dx
+\int_{\R^d}  \nabla [(p^\e)^{\frac{1}{\gamma}+1}] c_i^\e\nabla \phi dx.
\end{aligned}
\end{equation}
Since $p^\e$ and $c_i^\e$ are bounded, then the first term on the right-hand side is bounded due to \eqref{est6}. The boundedness of the second term is provided by \eqref{est2}. For the third and fifth term we use Young's inequality and argue with boundedness of $\nabla p^\e$ in $L^2_{t,x}$. The fourth term is obvious. Thus after letting $\e\to 0$ this error term will vanish. Finally we let  $\delta\to 0$  and  complete the proof. \qed

%%%%%%%%%%%%
\section{Conclusion and perspectives}
%-----------------------------
%%%%%%%%%%%%

We have proposed a strategy to prove existence of weak solutions for a two species model of tumor invasion. It relies on the extension of  the Aronson-Benilan regularizing effect for porous media equations which provides estimates of the Laplacian of the pressure. The most important limitation so far is a combined condition on the two bulk growth terms and it is an open question to remove it.  A route in this direction  could be to use the energy  type estimate given in Theroem~\ref{th:Aenergy} in the appendix. 
\\

A question which we do not handle here is the strong compactness on the $n^\e_i$ in the stability result of the approximation process. The bounds on $\Delta p$ are too weak for the $L^1$  theory in  \cite{DPL89} and are boarder line to apply the compactness theorems in \cite{Ambrosio04, BoCr06} which require that $D^2 p $ is a bounded measure.
\\

The extension to more than two species, with the present strategy, requires combined conditions on the three growth terms which read, in the case of three species for instance, $c_1 F(p) +c_2 G(p) +c_3H(p) \leq Cp^{1/\gamma}$ whenever the nonnegative $c_i$ satisfy $c_1+c_2+c_3=1$.  Then, the analysis goes through without major changes.
 \\

There are other questions which arise in this area and that we leave open. One of them concerns the `"incompressible  limit' $\gamma \to \infty$ which has attracted much attention recently \cite{PQV, KTu,MR3740370,Kpo,KPS} because of its relation to congested traffic \cite{MRS1,MR2651987,MRS2,MRSV}.  Clearly the bounds provided here are not enough to investigate this question. However the one dimensional case is under investigation~\cite{BPPS} based upon arguments from \cite{CFSS}. Another question is about  different mobilities, see \cite{LLP,ByDr,byrneChaplain96}, where the parabolic aspects of the equation for $n=n_1+n_2$ do not apply. 

%%%%%%%%%
\appendix
%\section{Appendix}
%-----------------------------
%%%%%%%%%%%%

\section{Additional a priori bounds}

Another remarkable estimate can be obtained for solutions of the system~ \eqref{eq1}--\eqref{eq2}. We give it here for the sake of completeness. It can be interpreted as some kind of energy because the kinetic energy is given by $E_K= n \f{|v|^2}{2}= p^{1/\gamma} \f{|\nabla p|^2}{2}$.

\begin{theorem}[Energy type a priori estimates]
With the assumptions \eqref{as1}--\eqref{as4},  the following estimates hold true with constants $C(T)$ which only depend on the bounds in the above assumptions. For $ \al_* = \f 2 \gamma$, we control 
\beq 
\dis \int_0^T \int_{\R^d}   \left[ \dv(p^\f{\al_*+1}{2}  \nabla p) - p^\f{\al_*+1}{2}\f{|\nabla p|^2}{2p} \right]^2 \leq C(T),
 \label{est3}
\eeq
\beq 
\dis \int_{\R^d} p^{\al_*} \;  |\nabla p (t) |^2  \,dx \leq C(T) \qquad \forall t\in (0,T).
\label{est4}
\eeq
\label{th:Aenergy}
\end{theorem}

\noindent
{\bf Proof}. 
These two estimates, \eqref{est3} and \eqref{est4} come together and require some elaborate computations.  We write
$$
\p_t \nabla p = \nabla [ |\nabla p|^2+  \gamma p \Delta p + \gamma p R],
$$
$$
\p_t \f{|\nabla p|^2}{2} = \nabla p .  \nabla [ |\nabla p|^2+  \gamma p \Delta p + \gamma p R],
$$

$$
\p_t p^\al \f{|\nabla p|^2}{2} =p^\al  \nabla p .  \nabla [ |\nabla p|^2+  \gamma p \Delta p + \gamma p R] + \al p^{\al -1}  \f{|\nabla p|^2}{2} [ |\nabla p|^2+  \gamma p \Delta p + \gamma p R]  .
$$
Therefore, we find
$$
\f {d}{dt} \int_{\R^d} p^\al \f{|\nabla p|^2}{2} = \int_{\R^d}[- p^\al  \Delta p- \al p^{\al-1}  |\nabla p|^2 +\al p^{\al -1}  \f{|\nabla p|^2}{2}] \;  [ |\nabla p|^2+  \gamma p \Delta p + \gamma p R] 
$$
$$
\f {d}{dt} \int_{\R^d} p^\al \f{|\nabla p|^2}{2} = \int_{\R^d}[- p^\al  \Delta p- \al p^{\al -1}  \f{|\nabla p|^2}{2}] \;  [ |\nabla p|^2+  \gamma p \Delta p + \gamma p R] 
$$
but $ p^\al  \Delta p$ is not a good quantity. So the right-hand side has to be rewritten (divide it by $\gamma$)
$$
- [ p^\f{\al+1}{2}  \Delta p + \al p^\f{\al+1}{2}   \f{|\nabla p|^2}{2p}] \;  [p^\f{\al+1}{2} \f{ |\nabla p|^2}{p\gamma}+   p^\f{\al+1}{2}   \Delta p +  p^\f{\al+1}{2}  R]  =
$$
$$
- [ \dv(p^\f{\al+1}{2}  \nabla p) - \f{\al+1}{2} p^\f{\al+1}{2}\f{|\nabla p|^2}{p}  + \al p^\f{\al+1}{2}   \f{|\nabla p|^2}{2p}] \; 
 [p^\f{\al+1}{2} \f{ |\nabla p|^2}{p\gamma}-   \f{\al+1}{2} p^\f{\al+1}{2}\f{|\nabla p|^2}{p}+   \dv(p^\f{\al+1}{2}  \nabla p)  +  p^\f{\al+1}{2}  R] 
$$
$$
= - [ \dv(p^\f{\al+1}{2}  \nabla p) -  p^\f{\al+1}{2}\f{|\nabla p|^2}{2p} ] \;  
[   \dv(p^\f{\al+1}{2}  \nabla p)  - p^\f{\al+1}{2} \f{ |\nabla p|^2}{2p} ( \al +1-\f 2 \gamma ) +  p^\f{\al+1}{2}  R] .
$$
To create a negative square, we use  the special value of $\al$ given by
$$
 \al_* = \f 2 \gamma
$$
and the right-hand side is controlled as $- [ \dv(p^\f{\al+1}{2}  \nabla p) -  p^\f{\al+1}{2}\f{|\nabla p|^2}{2p} ]^2 + C \big| \dv(p^\f{\al+1}{2}  \nabla p) -  p^\f{\al+1}{2}\f{|\nabla p|^2}{2p} \big|$.

Therefore, we obtain the inequalities  \eqref{est3} and \eqref{est4}.
\qed

%
%%%%%%%%%%%%%%%%%%%%%%%%%%%%%%%%%%%
%
%%%%%% BIBLIO %%%%%%%%%%%%%%%%%%%%%%
%
%%%%%%%%%%%%%%%%%%%%%%%%%%%%%%%%%%%%
%\pagestyle{myheadings}

\bibliographystyle{abbrv}
\bibliography{cancer}

\begin{thebibliography}{10}

\bibitem{ADiB}
H.~W. Alt and E.~DiBenedetto.
\newblock Nonsteady flow of water and oil through inhomogeneous porous media.
\newblock {\em Ann. Scuola Norm. Sup. Pisa Cl. Sci. (4)}, 12(3):335--392, 1985.

\bibitem{ALuVi}
H.~W. Alt, S.~Luckhaus, and A.~Visintin.
\newblock On nonstationary flow through porous media.
\newblock {\em Ann. Mat. Pura Appl. (4)}, 136:303--316, 1984.

\bibitem{Ambrosio04}
L.~Ambrosio.
\newblock Transport equation and {C}auchy problem for {$BV$} vector fields.
\newblock {\em Invent. Math.}, 158(2):227--260, 2004.

\bibitem{AB}
D.~G. Aronson and P.~B{\'e}nilan.
\newblock R{\'e}gularit{\'e} des solutions de l'{\'e}quation des milieux poreux
  dans {${\bf R}^{N}$}.
\newblock {\em C. R. Acad. Sci. Paris S{\'e}r. A-B}, 288(2):A103--A105, 1979.

\bibitem{BGH}
M.~Bertsch, M.~E. Gurtin, and D.~Hilhorst.
\newblock On interacting populations that disperse to avoid crowding: the case
  of equal dispersal velocities.
\newblock {\em Nonlinear Anal.}, 11(4):493--499, 1987.

\bibitem{BHIM}
M.~Bertsch, D.~Hilhorst, H.~Izuhara, and M.~Mimura.
\newblock A nonlinear parabolic-hyperbolic system for contact inhibition of
  cell-growth.
\newblock {\em Differ. Equ. Appl.}, 4(1):137--157, 2012.

\bibitem{BoCr06}
F.~Bouchut and G.~Crippa.
\newblock Uniqueness, renormalization, and smooth approximations for linear
  transport equations.
\newblock {\em SIAM J. Math. Anal.}, 38(4):1316--1328, 2006.

\bibitem{MR2651987}
L.~Brasco, G.~Carlier, and F.~Santambrogio.
\newblock Congested traffic dynamics, weak flows and very degenerate elliptic
  equations [corrected version of mr2584740].
\newblock {\em J. Math. Pures Appl. (9)}, 93(6):652--671, 2010.

\bibitem{BCGR}
D.~Bresch, T.~Colin, E.~Grenier, B.~Ribba, and O.~Saut.
\newblock Computational modeling of solid tumor growth: the avascular stage.
\newblock {\em SIAM J. Sci. Comput.}, 32(4):2321--2344, 2010.

\bibitem{BPPS}
F.~Bubba, C.~Pouchol, B.~Perthame, and M.~Schmidtchen.
\newblock Incompressible limit for a two species model of tissue growth in one
  space dimension.
\newblock Work in preparation.

\bibitem{byrneChaplain96}
H.~M. Byrne and M.~Chaplain.
\newblock Growth of necrotic tumors in the presence and absence of inhibitors.
\newblock {\em Mathematical biosciences}, 135(2):187--216, 1996.

\bibitem{ByDr}
H.~M. Byrne and D.~Drasdo.
\newblock Individual-based and continuum models of growing cell populations: a
  comparison.
\newblock {\em Math. Med. Biol.}, 58(4-5):657--687, 2003.

\bibitem{BKMP2}
H.~M. Byrne, J.~R. King, D.~L.~S. McElwain, and L.~Preziosi.
\newblock A two-phase model of solid tumour growth.
\newblock {\em Appl. Math. Lett.}, 16(4):567--573, 2003.

\bibitem{BP}
H.~M. Byrne and L.~Preziosi.
\newblock Modelling solid tumour growth using the theory of mixtures.
\newblock {\em Math. Med. Biol.}, 20(4):341--366, 2003.

\bibitem{CGM}
C.~Canc\`es, T.~O. Gallou\"et, and L.~Monsaingeon.
\newblock Incompressible immiscible multiphase flows in porous media: a
  variational approach.
\newblock {\em Anal. PDE}, 10(8):1845--1876, 2017.

\bibitem{CFSS}
J.~A. Carrillo, S.~Fagioli, F.~Santambrogio, and M.~Schmidtchen.
\newblock Splitting schemes \& segregation in reaction-(cross-)diffusion
  systems.
\newblock {\em arxiv:1711.05434}, 2017.

\bibitem{MR3740370}
K.~Craig, I.~Kim, and Y.~Yao.
\newblock Congested aggregation via {N}ewtonian interaction.
\newblock {\em Arch. Ration. Mech. Anal.}, 227(1):1--67, 2018.

\bibitem{DPL89}
R.~J. DiPerna and P.-L. Lions.
\newblock Ordinary differential equations, transport theory and {S}obolev
  spaces.
\newblock {\em Invent. Math.}, 98(3):511--547, 1989.

\bibitem{GorielyBenAmar2005}
A.~Goriely and M.~B. Amar.
\newblock Differential growth and instability in elastic shells.
\newblock {\em Phys. Rev. Lett.}, 94:198103--1 -- 198103--4, 2005.

\bibitem{Kpo}
I.~Kim and N.~Po{\v z\'a}r.
\newblock Porous medium equation to {H}ele-{S}haw flow with general initial
  density.
\newblock {\em Trans. Amer. Math. Soc.}, 370(2):873--909, 2018.

\bibitem{KTu}
I.~Kim and O.~Turanova.
\newblock Uniform convergence for the incompressible limit of a tumor growth
  model.
\newblock {\em Ann. Inst. H. Poincar\'e Anal. Non Lin\'eaire},
  35(5):1321--1354, 2018.

\bibitem{KPS}
I.~C. Kim, B.~Perthame, and P.~E. Souganidis.
\newblock Free boundary problems for tumor growth: a viscosity solutions
  approach.
\newblock {\em Nonlinear Anal.}, 138:207--228, 2016.

\bibitem{LLP}
T.~Lorenzi, A.~Lorz, and B.~Perthame.
\newblock On interfaces between cell populations with different mobilities.
\newblock {\em Kinet. Relat. Models}, 10(1):299--311, 2017.

\bibitem{MRS1}
B.~Maury, A.~Roudneff-Chupin, and F.~Santambrogio.
\newblock A macroscopic crowd motion model of gradient flow type.
\newblock {\em Math. Models Methods Appl. Sci.}, 20(10):1787--1821, 2010.

\bibitem{MRS2}
B.~Maury, A.~Roudneff-Chupin, and F.~Santambrogio.
\newblock Congestion-driven dendritic growth.
\newblock {\em Discrete Contin. Dyn. Syst.}, 34(4):1575--1604, 2014.

\bibitem{MRSV}
B.~Maury, A.~Roudneff-Chupin, F.~Santambrogio, and J.~Venel.
\newblock Handling congestion in crowd motion modeling.
\newblock {\em Netw. Heterog. Media}, 6(3):485--519, 2011.

\bibitem{PQTV}
B.~Perthame, F.~Quir{\'o}s, M.~Tang, and N.~Vauchelet.
\newblock Derivation of a {H}ele-{S}haw type system from a cell model with
  active motion.
\newblock {\em Interfaces Free Bound.}, 16(4):489--508, 2014.

\bibitem{PQV}
B.~Perthame, F.~Quir{\'o}s, and J.-L. V{\'a}zquez.
\newblock The {H}ele-{S}haw asymptotics for mechanical models of tumor growth.
\newblock {\em Arch. Ration. Mech. Anal.}, 212(1):93--127, 2014.

\bibitem{PT}
L.~Preziosi and A.~Tosin.
\newblock Multiphase modelling of tumour growth and extracellular matrix
  interaction: mathematical tools and applications.
\newblock {\em J. Math. Biol.}, 58(4-5):625--656, 2009.

\bibitem{RJPJ}
J.~Ranft, M.~Basana, J.~Elgeti, J.-F. Joanny, J.~Prost, and F.~J\"ulicher.
\newblock Fluidization of tissues by cell division and apoptosis.
\newblock {\em Natl. Acad. Sci. USA}, 49:657--687, 2010.

\bibitem{RSCB}
B.~Ribba, O.~Saut, T.~Colin, D.~Bresch, E.~Grenier, and J.~P. Boissel.
\newblock A multiscale mathematical model of avascular tumor growth to
  investigate the therapeutic benefit of anti-invasive agents.
\newblock {\em J. Theoret. Biol.}, 243(4):532--541, 2006.

\bibitem{SCh}
J.~A. Sherratt and M.~A.~J. Chaplain.
\newblock A new mathematical model for avascular tumour growth.
\newblock {\em J. Math. Biol.}, 43(4):291--312, 2001.

\bibitem{Vbook}
J.-L. V{\'a}zquez.
\newblock {\em The porous medium equation. Mathematical theory}.
\newblock Oxford Mathematical Monographs. The Clarendon Press, Oxford
  University Press, Oxford, 2007.

\end{thebibliography}

\end{document}